# Continuous rotation invariant valuations
# on convex sets

By S. Alesker*

## 1. Introduction

The notion of valuation on convex sets can be considered as a generalization of the notion of measure, which is defined only on the class of convex compact sets. It is well-known that there are important and interesting examples of valuations on convex sets, which are not measures in the usual sense as, for example, the mixed volumes. Basic definitions and some classical examples are discussed in Section 2 of this paper. For more detailed information we refer to the surveys [Mc-Sch] and [Mc3]. Throughout this paper all the valuations are assumed to be continuous with respect to the Hausdorff metric.

Note that the theory of valuations which are invariant or covariant with respect to translations belongs to the classical part of convex geometry. There exists an explicit description of translation invariant continuous valuations on $\mathbb{R}^1$ and $\mathbb{R}^2$ due to Hadwiger [H1] (the case of $\mathbb{R}^2$ is nontrivial). Continuous rigid motion invariant valuations on $\mathbb{R}^d$ are completely classified by the remarkable Hadwiger theorem as linear combinations of the *quermassintegrals* (cf. [H2] or for a simpler proof [K]).

There are two natural ways to generalize Hadwiger's theorem: the first one is to describe continuous translation invariant valuations without any assumption on rotations; the second one is to characterize continuous rotation (i.e. either O($d$)- or SO($d$)-) invariant valuations without any assumption on translations (here O($d$) denotes the full orthogonal group and SO($d$) denotes the special orthogonal group). The first problem is of interest to classical convexity and translative integral geometry. As we have said, it was solved by Hadwiger for the line and the 2-dimensional plane. There is a conjecture due to P. McMullen [Mc2], which states that every continuous translation invariant valuation can be approximated (in some sense) by linear combinations of mixed volumes (note that in the 3-dimensional space this conjecture is known to be true and it follows from several other general results, which we do not discuss here).

*Partially supported by a Minkowski Center grant and BSF grant.

The main goal of this paper is to solve the second problem, namely to present a characterization of continuous O($d$)- (resp. SO($d$)-) invariant valuations.

Originally the second problem was motivated by questions arising in the asymptotic theory of normed spaces, where the property of invariance with respect to rotations is more natural than that of invariance with respect to translations. For example the following expression (which is a valuation in $K$) is of great interest in the asymptotic theory

$$\varphi(K) = \int\limits_K |x|^2 dx \ ,$$

where $K$ is a convex compact set, and $|\cdot|$ is the Euclidean norm. For detailed discussion we refer to [M-P]; see also [Bo].

The space of all continuous rotation invariant valuations is infinite-dimensional. To describe it, we consider a smaller subspace of *polynomial* continuous rotation invariant valuations (see Definition 2.2 below), which turns out to be everywhere dense and which has a natural filtration with respect to the degree of polynomiality.

The class of polynomial valuations was introduced by Khovanskii and Pukhlikov [Kh-P1] for polytopes. They developed the combinatorial theory of these valuations, which was applied in the subsequent paper [Kh-P2] to obtain a Riemann-Roch type theorem for integrals and sums of quasipolynomials over polytopes.

Let us denote by $\mathcal{K}^d$ the family of convex compact subsets of $\mathbb{R}^d$. Equipped with the Hausdorff metric, $\mathcal{K}^d$ is a locally compact space. Our first main result is:

THEOREM A. *Every continuous* SO($d$)- (*resp.* O($d$)-) *invariant valuation can be approximated uniformly on compact subsets in* $\mathcal{K}^d$ *by continuous polynomial* SO($d$)- (*resp.* O($d$)-) *invariant valuations.*

Thus the problem of describing continuous rotation invariant valuations is reduced to a more natural one of describing *polynomial* continuous rotation invariant valuations. Our second main result states that such valuations can be described explicitly by presenting a complete list of them. The linear space of polynomial continuous O($d$)- (resp. SO($d$)-) invariant valuations has the natural increasing filtration with respect to the degree of polynomiality. In particular it is shown that the space of valuations, which are polynomial up to a given degree, is finite dimensional.

In order to state precisely our second main theorem, we will need the notion of the generalized curvature measure of a convex set $K$, for the definition of which we refer to [Sch1]. However, this is not strictly necessary for the statement of the theorem, and the reader who feels uncomfortable with this terminology can find an equivalent form of the main result in Theorems 4.7 and 4.4 below (but then the formulation becomes longer). So let us denote by $\Theta_j(K; \cdot)$ the $j^{\text{th}}$ generalized curvature measure of $K$, which is defined on $\mathbb{R}^d \times S^{d-1}$, where $S^{d-1}$ in the unit sphere in $\mathbb{R}^d$. Then we have:

THEOREM B.    (i) *Let $\varphi$ be a continuous polynomial valuation, which is* SO($d$)-*invariant if $d \geq 3$ and* O($d$)-*invariant if $d = 2$. Then there exist polynomials $p_0, \ldots, p_{d-1}$ in two variables such that*

$$(1.1) \qquad \varphi(K) = \sum_{j=0}^{d-1} \int_{\mathbb{R}^d \times S^{d-1}} p_j(|s|^2, \langle s, n \rangle) d\Theta_j(K; s, n)$$

*for every $K \in \mathcal{K}^d$, where $\Theta_j(K; \cdot)$ is the $j^{\text{th}}$ generalized curvature measure of $K$, $|s|$ is the Euclidean norm of $s \in \mathbb{R}^d$, and $n \in S^{d-1}$. Moreover, any expression of the form* (1.1) *is a continuous polynomial* O($d$)-*invariant valuation.*

(ii) *Let $\varphi$ be a continuous polynomial* SO(2)-*invariant valuation on $\mathcal{K}^2$. Then there exist polynomials $q_0, q_1$ in two variables such that*

$$(1.2) \qquad \varphi(K) = \sum_{j=0}^{1} \int_{\mathbb{R}^2 \times S^1} q_j(\langle s, n \rangle, \langle s, n' \rangle) d\Theta_j(K; s, n) ,$$

*where $n'$ denotes the rotation of the vector $n$ by $\frac{\pi}{2}$ counterclockwise. Moreover, any expression of the form* (1.2) *is a continuous* SO(2)-*invariant polynomial valuation.*

Thus Theorems A and B give a complete description of all continuous rotation invariant valuations in $\mathbb{R}^d$. Note also an immediate corollary of Theorem A and Theorem B(i): in dimension $d \geq 3$ every continuous SO($d$)-invariant valuation is O($d$)-invariant (but this is not true if $d \leq 2$). We do not know any direct explanation of this corollary.

The paper is organized as follows. Section 2 contains necessary definitions, examples and known results used in the paper.

In Section 3 we present a description of valuations on the line (which is in fact rather trivial).

Section 4 contains the proof of the main Theorems A and B.

In Section 5 we give some applications of the main results to integral-geometric formulas.

In Section 6 we discuss some inequalities related to concrete polynomial valuations. Thus Theorem 6.1 says that the polynomial $\int_{K + \varepsilon B} |s|^{2q} dm(s)$ has nonnegative coefficients in $\varepsilon \geq 0$, where $K$ is a convex compact set containing the origin, $B$ is the Euclidean ball, and $q$ is a nonnegative integer.

In Section 7 we state several natural questions.

*Remark.* After the preprint of this paper was distributed we received from Prof. P. McMullen a preprint of his work [Mc4], where a more general class of valuations was introduced (isometry covariant valuations) . Some concrete examples of valuations and relations between them were studied and there was formulated a conjecture on characterization of such valuations. The methods of our paper turned out to be useful in solving this conjecture (see [A2]).

*Acknowledgements.* We are grateful to Professor Vitali Milman for his guidance in this work. We would also like to thank Professor J. Bernstein and Dr. A. Litvak for useful discussions, and Professor P. McMullen for important remarks.

## 2. Preliminaries

In this section we present some notation, definitions and facts used in the paper.

Let $\mathcal{K}^d$ denote the family of all compact convex subsets of $\mathbb{R}^d$. Let $L$ be a finite dimensional vector space over $\mathbb{R}$ or $\mathbb{C}$.

*Definition* 2.1. A function $\varphi : \mathcal{K}^d \to L$ is called a *valuation*, if $\varphi(K_1 \cup K_2) + \varphi(K_1 \cap K_2) = \varphi(K_1) + \varphi(K_2)$ for all $K_1, K_2 \in \mathcal{K}^d$ such that $K_1 \cup K_2 \in \mathcal{K}^d$. If $\varphi$ is continuous with respect to the Hausdorff metric on $\mathcal{K}^d$, we call it a *continuous valuation*; only such valuations will be considered here.

*Definition* 2.2. The valuation $\varphi : \mathbb{R}^d \to L$ is called *polynomial* of degree at most $\ell$ if $\varphi(K + x)$ is a polynomial in $x$ of degree at most $\ell$ for each $K \in \mathcal{K}^d$. Valuations of degree 0 correspond to the translation invariant valuations, and those of degree 1 to translation covariant ones.

The following theorem due to Khovanskii and Pukhlikov [Kh-P1] (actually, a special case) generalizes to the polynomial case the previous result of McMullen [Mc1] obtained for translation invariant and translation covariant valuations.

THEOREM 2.3. *Let* $\varphi : \mathcal{K}^d \to L$ *be a continuous valuation, which is polynomial of degree at most* $\ell$. *Then, for every* $K_1, \ldots, K_s \in \mathcal{K}^d$, $\varphi\left(\sum_{j=1}^{s} \lambda_j K_j\right)$ *is a polynomial in* $\lambda_j \geq 0$ *of degree at most* $d + \ell$, *where* $\sum_{j=1}^{s} \lambda_j K_j$ *denotes the Minkowski sum of the sets* $\lambda_j K_j$.

For the proof of Theorem 2.3, see [Kh-P1] or [A1].

Now let us recall some well-known results on translation invariant valuations, which will be used in the sequel.

THEOREM 2.4.    *Let $\varphi : \mathcal{K}^d \to \mathbb{R}$ be a continuous translation invariant valuation.*

(a) *If $d = 1$, then $\varphi$ has the form $\varphi(K) = a + b|K|$, where $|K|$ is the length of $K \in \mathcal{K}^1$ (i.e. the Lebesgue measure of $K$), and $a, b$ are uniquely defined constants.*

(b) *(Hadwiger [H1]) If $d = 2$, then $\varphi$ has the form*

$$\varphi(K) = a + b \ \mathrm{vol}_2 K + \int\limits_{S^1} f(\omega) dS_1(K, \omega) \ ,$$

*where $a, b$ are constants, $f : S^1 \to \mathbb{R}$ is a continuous function on the unit circle, and $S_1(K, \cdot)$ is the surface area measure of $K$.*

(c) *(Hadwiger [H2]) If, in addition, $\varphi$ is $\mathrm{SO}(d)$-invariant, then $\varphi$ has the form*

$$\varphi(K) = \sum_{j=1}^{d} c_j W_j(K) \ ,$$

*where $W_j(K)$ is the $j^{\mathrm{th}}$ quermassintegral, and the $c_j$ are fixed, uniquely defined constants.*

For the definition of the surface area measure and the quermassintegrals we refer to [Sch1]. Obviously, Theorem 2.4 generalizes immediately to $L$-valued valuations.

The following result is an easy consequence of the translation invariant (McMullen's) version of Theorem 2.3.

THEOREM 2.5 ([Mc1]).    *Let $\varphi : \mathcal{K}^d \to L$ be a continuous translation invariant valuation. Then $\varphi$ can be uniquely represented as a sum*

$$\varphi = \sum_{j=0}^{d} \varphi_j \ ,$$

*where $\{\varphi_j\}$ are translation invariant continuous valuations, homogeneous of degree $j$, so that for every $K \in \mathcal{K}^d$ and every $\lambda \geq 0$,*

$$\varphi_j(\lambda K) = \lambda^j \varphi_j(K) \ .$$

THEOREM 2.6.    *Let $\varphi : \mathcal{K}^d \to \mathbb{R}$ be a continuous translation invariant valuation, homogeneous of degree $j$. Then*

(a) *(trivial) $\varphi_0$ is just a constant;*

(b) *([H2]) $\varphi_d$ is a multiple of the standard volume; i.e. $\varphi_d(K) = a \cdot \mathrm{vol}_d K$;*

(c)  ([Mc2])

$$\varphi_{d-1}(K) = \int\limits_{S^{d-1}} f(\omega) dS_{d-1}(K, \omega) \ ,$$

*where $f : S^{d-1} \to \mathbb{R}$ is a continuous function, and $S_{d-1}(K, \cdot)$ is the surface area measure of $K$.*

*Remark* 1.  It is well-known that the function $f$ in 2.6(c) and 2.4(b) can be chosen to be orthogonal to every linear functional on $\mathbb{R}^d$ with respect to the standard Lebesgue (Haar) measure on $S^{d-1}$. Under this assumption $f$ is unique (this follows from Minkowski's existence theorem; cf. e.g. [Mc2, Th. 3]).

*Remark* 2.  Theorem 2.4(a) and (b) immediately follow from Theorems 2.5 and 2.6.

The next theorem was recently established by Schneider [Sch2], but a particular case of the even valuations was considered by Klain [K].

THEOREM 2.7.    *Let $\varphi : \mathcal{K}^d \to \mathbb{R}$ be a continuous translation invariant valuation, which is simple, i.e. $\varphi(K) = 0$ whenever $\dim K < d$. Then $\varphi$ has the form*

$$\varphi(K) = a \operatorname{vol}_d K + \int\limits_{S^{d-1}} f(\omega) dS_{d-1}(K, \omega) \ ,$$

*where $f$ is a continuous odd function on the unit sphere.*

Again, $f$ has the same uniqueness properties as in Remark 1 above.

Now we will give some examples of rotation invariant polynomial valuations. Fix a nonnegative integer $m$ and consider $\varphi : \mathcal{K}^d \to \mathbb{R}$ given by

$$\varphi(K) := \int\limits_{K} |x|^{2m} dx \ ,$$

where $|x|$ denotes the Euclidean norm of $x \in \mathbb{R}^d$. Then obviously $\varphi$ is a continuous O$(d)$-invariant valuation, polynomial of degree $2m$. It is well known (e.g. [Sch1, p. 173]) that if $\varphi$ is a valuation and $A \in \mathcal{K}^d$ is fixed, then $\psi(K) := \varphi(K + A)$ is also a valuation. Thus for every $\varepsilon \geq 0$ we may consider the valuation $\varphi(K + \varepsilon B)$, where $B$ is the Euclidean ball in $\mathbb{R}^d$. Clearly, this is also a continuous O$(d)$-invariant valuation, polynomial of degree $2m$. But by Theorem 2.3 this is a polynomial in $\varepsilon$ (for a fixed $K$) of degree $2m + d$, whose coefficients are also continuous O$(d)$-invariant polynomial valuations. Thus $\left( \frac{d^j}{d\varepsilon^j} \right)\Big|_{\varepsilon=0} \varphi(K + \varepsilon B)$, $0 \leq j \leq 2m + d$, gives us more examples of such valuations.

Later on we will show that there are other valuations which cannot be expressed as linear combinations of valuations of the above type.

## 3. Polynomial valuations on the line $\mathbb{R}^1$

PROPOSITION 3.1.    *Every continuous valuation $\varphi : \mathcal{K}^1 \to \mathbb{C}$ has the form*

$$\varphi([a,b]) = P(a) + Q(b),$$

*for every segment $[a,b] \subset \mathbb{R}^1$, where $P, Q$ are continuous functions on $\mathbb{R}^1$. Moreover, if $\varphi$ is a polynomial valuation of degree at most $\ell$, then $P, Q$ can be chosen to be polynomials of degree at most $\ell + 1$.*

*Proof.* Let us prove this for polynomial valuations. By definition, $\varphi(\{x\}) = \varphi(\{0\} + x)$ is a polynomial $T(x)$ in $x$ of degree at most $\ell$. Therefore the valuation $\psi([a,b]) := \varphi([a,b]) - T(a)$ vanishes on points. By Theorem 2.3, $\psi([0,x])$, $x \geq 0$, is a polynomial of degree at most $\ell$. Denote it by $S(x)$, and then obviously $\psi([a,b]) = S(b) - S(a)$ for every $a \leq b$. Thus $\varphi([a,b]) = S(b) - (S - T)(a)$.    □

Note that the functions $P, Q$ in Proposition 3.1 are defined uniquely up to the same constant.

## 4. Main results: Rotation invariant polynomial valuations in dimension greater than 1

Since a linear combination of valuations is again a valuation, we will denote by $\Omega_{d,\ell}$(resp. $\Omega'_{d,\ell}$) the linear space of continuous SO($d$)- (resp. O($d$)-) invariant valuations on $\mathcal{K}^d$, which are polynomial of degree at most $\ell$. Clearly,

$$\Omega_{d,\ell} \supset \Omega'_{d,\ell} \ ,$$
$$\Omega_{d,0} \subset \Omega_{d,1} \subset \cdots \subset \Omega_{d,\ell} \subset \cdots$$

and the similar sequence of inclusions holds for $\Omega'_{d,\ell}$.

The first result of this section is:

THEOREM 4.1.    $\Omega_{d,1} = \Omega_{d,0}(= \Omega'_{d,1})$ *if $d \geq 2$.*

Before proving this result, we observe that if $\varphi$ is a polynomial valuation of degree $\ell$, then for every $K \in \mathcal{K}^d$,

$$\varphi(K + x) = P_K^\ell(x) + P_K^{\ell-1}(x) + \cdots + P_K^0(x) \ ,$$

where $P_K^j(x)$ is a homogeneous polynomial of degree $j$ with coefficients depending on $K$. Then the $P_K^j$ have the following properties.

(i)  $P_K^j$ is a continuous (polynomial-valued) valuation in $K$ (i.e. $P_{K_1 \cup K_2}^j + P_{K_1 \cap K_2}^j \equiv P_{K_1}^j + P_{K_2}^j$ whenever $K_1, K_2, K_1 \cup K_2 \in \mathcal{K}^d$);

(ii) $P_K^\ell$ is a translation invariant valuation;

(iii) If $\varphi$ is SO($d$)- (resp. O($d$)-) invariant, then $P_K^\ell$ is SO($d$)- (resp. O($d$)-) equivariant; i.e., for every $U \in$ SO($d$) (resp. O($d$)) and every $K \in \mathcal{K}^d$,

$$(4.1) \qquad\qquad P_{UK}^\ell \equiv \pi(U)(P_K^\ell) \ ,$$

where $\pi(U)$ denotes the standard quasi-regular representation of SO($d$) (resp. O($d$)) in the space of homogeneous polynomials in $d$ variables of degree $\ell$ acting as $\pi(U)P_K^\ell(x) = P_K^\ell(U^{-1}x)$ (cf. [V]);

Let us check, for example, (ii). Fix $y \in \mathbb{R}^d$, and $K \in \mathcal{K}^d$. Then $\varphi((K+y) + x) = P_{K+y}^\ell(x)+$ (lower order terms), where the expression "lower order terms" means the sum of monomials in $x$ of degree strictly less than $l$. However, the left-hand side equals

$$\varphi(K + (x+y)) = P_K^\ell(x+y) + \text{ (lower order terms)}$$
$$= P_K^\ell(x) + \text{ (lower order terms) .}$$

Comparing the right-hand sides of these expressions, we get (ii).

We denote by $T_{d,\ell}$ the finite-dimensional space of homogeneous polynomials in $d$ variables of degree $\ell$, and by $\Gamma_{d,\ell}$ (resp. $\Gamma'_{d,\ell}$) the linear space of $T_{d,\ell}$-valued valuations satisfying properties (i)–(iii). $\Gamma_{d,\ell}$ corresponds to the case of SO($d$) in (iii), and $\Gamma'_{d,\ell}$ to O($d$). Clearly, $\Gamma'_{d,\ell} \subset \Gamma_{d,\ell}$.

The correspondence $\varphi \mapsto P_K^\ell$ defines a linear map $D : \Omega_{d,\ell} \longrightarrow \Gamma_{d,\ell}$ (resp. $D : \Omega'_{d,\ell} \longrightarrow \Gamma'_{d,\ell}$). Obviously, $\mathrm{Ker}\, D = \Omega_{d,\ell-1}$ (resp. $\Omega'_{d,\ell-1}$).

Theorem 4.1 immediately follows from:

PROPOSITION 4.2.     $\Gamma_{d,1} = 0$ for $d \geq 2$.

*Proof.* We use induction in $d$. Clearly, $T_{d,1}$ is isomorphic to $\mathbb{C}^d = \mathbb{R}^d \otimes \mathbb{C}$ as a representation of SO($d$). First, let $d = 2$. Fix $\Phi \in \Gamma_{2,1}$. By Theorem 2.4 (b),

$$\Phi(K) = A + B \cdot \mathrm{vol}_2\, K + \int_{S^1} F(\omega)dS_1(K,\omega) \ ,$$

where $A, B \in T_{2,1}$ ($\simeq \mathbb{C}^2$), and $F : S^1 \to T_{2,1}$ is a continuous function which is orthogonal to every linear functional. The uniqueness of such a representation and the rotation equivariance (4.1) imply that for every $U \in$ SO(2),

$$(4.2) \qquad\qquad\qquad UA = A,$$
$$(4.3) \qquad\qquad\qquad UB = B,$$
$$(4.4) \qquad\qquad\qquad UF(\omega) = F(U\omega) \ .$$

It follows from (4.2) and (4.3) that $A = B = 0$. Using (4.4) we define an intertwining operator $\widetilde{F} : T_{2,1}^* \to C(S^1)$ between the dual of the quasi-regular repre-sentation of $\mathrm{SO}(2)$ (in the dual space of $T_{2,1}$) and the quasi-regular representation in $C(S^1)$, as follows: for every $\xi \in T_{2,1}^*$ and $\omega \in S^1$ let $\widetilde{F}(\xi)(\omega) = \langle \xi, F(\omega) \rangle$. Let us denote by $\sigma_d^j$ the space of spherical harmonics in $d$ variables of degree $j$. If $d = 2$, then $\sigma_2^0$ is one-dimensional and for $j \geq 1$, $\sigma_2^j = \sigma_2^{'j} \oplus \sigma_2^{''j}$, where $\sigma_2^{'j}$ is spanned by $e^{ij\theta}$ and $\sigma_2^{''j}$ is spanned by $e^{-ij\theta}$, $\theta \in S^1$. Since all $\sigma_2^{'j}$ and $\sigma_2^{''k}$ are pairwise nonequivalent representations of $\mathrm{SO}(2)$ and since $C(S^1) = \sigma_2^0 \oplus \sigma_2^1 \oplus \cdots$, $\sigma_2^j = \sigma_2^{'j} \oplus \sigma_2^{''j}$, and $T_{2,1}^* = \sigma_2^1 = \sigma_2^{'1} \oplus \sigma_2^{''1}$, by Schur's lemma, we get that $\widetilde{F}(T_{2,1}) \subset \sigma_2^1$. Namely, for every $\xi \in T_{2,1}^*$,

$$\langle \xi, F(\omega) \rangle \in \sigma_2^1 \ ;$$

i.e., it is a restriction of a linear functional to the sphere $S^1$. But the assumption of orthogonality of $F$ to every linear functional implies that $F \equiv 0$. Thus $\Phi \equiv 0$ for $d = 2$.

Now let $d > 2$. Fix $\Phi \in \Gamma_{d,1}$ and an orthogonal decomposition $\mathbb{R}^d = \mathbb{R}^{d-1} \oplus \mathbb{R}^1$. If $K \subset \mathbb{R}^{d-1}$, $\Phi(K) = (\Phi_1(K), \Phi_2(K))$, where $\Phi_1$ is a projection of $\Phi$ onto $\mathbb{C}^{d-1} = \mathbb{R}^{d-1} \otimes \mathbb{C}$, and $\Phi_2$ is a projection of $\Phi$ onto $\mathbb{C}^1 = \mathbb{R}^1 \otimes \mathbb{C}$. Then clearly the restriction of $\Phi_1$ to $\mathcal{K}^{d-1}$ belongs to $\Gamma_{d-1,1}$, and the restriction of $\Phi_2$ to $\mathcal{K}^{d-1}$ is a $\mathbb{C}$-valued translation and $\mathrm{SO}(d-1)$-invariant valuation on $\mathcal{K}^{d-1}$. Thus, by the inductive assumption, $\Phi_1(K) = 0$ if $K \subset \mathbb{R}^{d-1}$, and the restriction of $\Phi_2$ to $\mathcal{K}^{d-1}$ satisfies Hadwiger's theorem 2.4(c). In particular, it is $\mathrm{O}(d-1)$-invariant.

Consider the following transformation $U \in \mathrm{SO}(d)$:

$$U(x_1, \ldots, x_{d-2}, x_{d-1}, x_d) = (x_1, \ldots, x_{d-2}, -x_{d-1}, -x_d) \ .$$

For every $K \subset \mathbb{R}^{d-1}$,

$$\Phi(UK) = (0, \Phi_2(UK)) = (0, \Phi_2(K)) \ ,$$
$$\Phi(UK) = U\Phi(K) = (0, -\Phi_2(K)) \ .$$

Thus $\Phi_2(K) = 0$; hence $\Phi(K) = 0$ for $K \subset \mathbb{R}^{d-1}$. By translation invariance and rotation equivariance, $\Phi$ vanishes on all $K$ such that $\dim K \leq d-1$. Hence by Theorem 2.7, it has the form

$$\Phi(K) = A \operatorname{vol}_d K + \int\limits_{S^{d-1}} F(\omega) dS_{d-1}(K, \omega) \ ,$$

where $F : S^{d-1} \to T_{d,1}$, and similarly to the 2-dimensional case (by pairwise nonequivalence of the $\sigma_d^j$) we deduce that $\Phi \equiv 0$. $\qquad\square$

Now we are going to give new examples of polynomial rotation invariant valuations. There is a difference between $\mathrm{SO}(d)$- and $\mathrm{O}(d)$-invariant valuations for $d = 2$ (for $d > 2$ we will see that there is no such difference).

Let $K \in \mathcal{K}^d$. For almost every point $s \in \partial K$, the unit outer normal $n(s)$ is defined uniquely. First consider the case $d = 2$. Denote by $n'(s)$ the rotation of $n(s)$ by the angle $\frac{\pi}{2}$ counterclockwise. Then define

$$(4.5) \qquad \psi_{p,q}(K) = \int_{\partial K} \langle s, n(s) \rangle^p \langle s, n'(s) \rangle^q d\sigma_K(s) \ ,$$

where $\sigma_K$ is the surface area measure on $\partial K$ and $p, q$ are fixed nonnegative integers. Note that if $K$ is a point, we set $\psi_{p,q}(K) = 0$.

PROPOSITION 4.3.    *The function* $\psi_{p,q}$ *is a continuous* SO(2)-*invariant polynomial valuation of degree of polynomiality* $\ell$, *where* $\ell = p + q$ *if* $p + q \neq 1$, *and* $\ell = 0$ *if* $p + q = 1$. *Moreover,* $\psi_{p,q}$ *is* O(2)-*invariant if and only if* $q$ *is even.*

*Proof.* Proof of the continuity of $\psi_{p,q}$ is standard. To see the valuation property, it is sufficient to check it in the following situation (see [G]): let $K \in \mathcal{K}^2$, $H$ be an affine hyperplane, and $H^+$ and $H^-$ be closed halfspaces into which $H$ divides $\mathbb{R}^2$. Then we have to verify

$$\psi_{p,q}(K) + \psi_{p,q}(K \cap H) = \psi_{p,q}(K \cap H^+) + \psi_{p,q}(K \cap H^-).$$

But this is immediate from the definition of $\psi_{p,q}$. Let us check polynomiality. Fix $K \in \mathcal{K}^2$. Then

$$\begin{aligned}
\psi_{p,q}(K + x) &= \int_{\partial K + x} \langle s, n(s) \rangle^p \langle s, n'(s) \rangle^q d\sigma_{K+x}(s) \\
&= \int_{\partial K} \langle s + x, n(s) \rangle^p \langle s + x, n'(s) \rangle^q d\sigma_K(s) \\
&= \int_{\partial K} \langle x, n(s) \rangle^p \langle x, n'(s) \rangle^q d\sigma_K(s) + \text{ lower order terms }.
\end{aligned}$$

Note that if $p + q = 1$ then the leading term above vanishes identically. However this does not happen if $p + q \neq 1$. If we denote by $R : \mathbb{R}^2 \to \mathbb{R}^2$ the rotation by $\frac{\pi}{2}$ counterclockwise, then $n'(s) = Rn(s)$, and the last integral can be rewritten as

$$\int_{\partial K} \langle x, n(s) \rangle^p \langle R^* x, n(s) \rangle^q d\sigma_K(s) = \int_{S^1} \langle x, \omega \rangle^p \langle R^* x, \omega \rangle^q dS_1(K, \omega) \ .$$

Thus if this expression is not identically 0 (for all $x \in \mathbb{R}^2$ and $K \in \mathcal{K}^2$), then $\psi_{p,q}$ has degree $p + q$. But if this expression vanishes, then by Remark 1 after Theorem 2.6 $\langle x, \omega \rangle^p \langle R^* x, \omega \rangle^q \equiv \langle \nu(x), \omega \rangle$, where $\nu(x)$ is a vector depending on $x$. This implies the first part of the proposition. The second part is clear. $\qquad \square$

THEOREM 4.4.    *Every continuous* SO(2)*-invariant polynomial valuation on* $\mathbb{R}^2$ *is a linear combination of valuations of the form* (4.5) *and of the form*

$$\frac{d^2}{d\varepsilon^2}\Big|_{\varepsilon=0} \int_{K+\varepsilon B} |x|^{2m} dx \ ,$$

*where* $m$ *is a nonnegative integer.*

*Remark.* One can easily see that Theorem 4.4 is equivalent to Theorem B (ii) in the introduction.

*Proof.* Let $\varphi : \mathcal{K}^2 \to \mathbb{C}$ satisfy the conditions of the theorem. Then $\varphi(\{x\})$ is a polynomial in $x \in \mathbb{R}^2$ which is SO(2)-invariant. Hence it has the form $\sum_{j\geq 0} c_j |x|^{2j}$. Consider a new valuation

$$\psi(K) = \varphi(K) - \sum_j \frac{c_j}{2 \cdot \mathrm{vol}_2 B} \frac{d^2}{d\varepsilon^2}\Big|_{\varepsilon=0} \int_{K+\varepsilon B} |x|^{2j} dx \ .$$

Clearly, $\psi$ vanishes on points. We will show that $\psi$ is a linear combination of valuations of the form (4.5).

Assume that $\psi$ has degree $\ell$. Recall that there is a map $D : \Omega_{2,\ell} \to \Gamma_{2,\ell}$. Theorem 4.4 follows by induction in $\ell$ from the following:

LEMMA 4.5.    *The span of* $\{D(\psi_{p,q})\}$ *coincides with all the valuations from* $\Gamma_{2,\ell}$ *vanishing on points (here* $p,q$ *are such that* $p+q = \ell$ *or* $p+q = \ell+1$ *and* $\psi_{p,q}$ *is a polynomial valuation of degree* $\ell$*).*

*Proof.* The case $\ell = 1$ follows from Proposition 4.2. Let $\ell > 1$, and fix $\Phi \in \Gamma_{2,\ell}$ vanishing on points. By Theorem 2.4(b),

$$\Phi(K) = A + B \ \mathrm{vol}_2(K) + \int_{S^1} F(\omega) dS_1(K,\omega) \ ,$$

where $A, B \in T_{2,\ell}$, $F : S^1 \to T_{2,\ell}$ is a continuous function which satisfies condition (4.1) of equivariance and $\pi(U)B = B$ for all $U \in$ SO(2). Since $\Phi$ vanishes on points, $A = 0$.

*Case* 1.    First we show that the valuation

$$\Phi_1(K) = \int_{S^1} F(\omega) dS_1(K,\omega)$$

belongs to the span of $\{D(\psi_{p,q})\}_{p+q=\ell}$.

Let us introduce the complex structure on the plane $\mathbb{R}^2$ in the standard way so that $1 = (1,0)$, $i = \sqrt{-1} = (0,1)$. Let $z$ denote a point of $\mathbb{C} \simeq \mathbb{R}^2$. For the quasi-regular representation of SO(2) in $T_{2,\ell}$ we have a decomposition into 1-dimensional (irreducible) components:

$$T_{2,\ell} = (\langle z^\ell \rangle \oplus \langle \overline{z}^\ell \rangle) \oplus (\langle z^{\ell-2} \cdot |z|^2 \rangle \oplus \langle \overline{z}^{\ell-2} \cdot |z|^2 \rangle) \oplus \cdots .$$

Using this decomposition, we may assume that $F$ takes values in the 1-dimensional space $\langle z^k |z|^{2m} \rangle$ or $\langle \overline{z}^k |z|^{2m} \rangle$, where $k + 2m = \ell$. Consider, e.g., the first case. Thus $F : S^1 \to \langle z^k |z|^{2m} \rangle$, and by equivariance $F$ has the form:

$$F(e^{i\theta}) = \alpha \cdot (ze^{-i\theta})^k |z|^{2m} ,$$

where $\alpha$ is some constant.

In the proof of Proposition 4.3 we have seen that

$$(D\psi_{p,q})(K)(z) = \int_{S^1} \langle z, \omega \rangle^p \langle R^* z, \omega \rangle^q dS_1(K, \omega) ,$$

where $p + q = \ell$ and $R^*$ is a rotation by the angle $-\frac{\pi}{2}$, so $R^* z = -i \cdot z$. Note that if $z, \omega \in \mathbb{C} \simeq \mathbb{R}^2$, then the scalar product $\langle z, \omega \rangle = \mathrm{Re}(z\overline{\omega})$, and $\langle R^* z, \omega \rangle = \mathrm{Im}(z\overline{\omega})$.

If we set $\omega = e^{i\theta} \in S^1$, then in this notation we obtain

$$(4.6) \qquad (D\psi_{p,q})(K)(z) = \int_{S^1} (\mathrm{Re} z e^{-i\theta})^p (\mathrm{Im} z e^{-i\theta})^q dS_1(K, \theta) .$$

Now we see that

$$\begin{aligned}
F(e^{i\theta}) &= \alpha (ze^{-i\theta})^k |z|^{2m} \\
&= \alpha (\mathrm{Re} z e^{-i\theta} + i\mathrm{Im} z e^{-i\theta})^k (|\mathrm{Re} z e^{-i\theta}|^2 + |\mathrm{Im} z e^{-i\theta}|^2)^m
\end{aligned}$$

belongs to the linear span of functions under the integral in (4.6). This implies Case 1.

*Case 2.* Let $\Phi_2(K) = B \cdot \mathrm{vol}_2(K)$, where $B \in T_{2,\ell}$ is an SO(2)-invariant polynomial. If $\ell$ is odd, then $B \equiv 0$ and there is nothing to prove. If $\ell$ is even, then $B$ has the form $B(z) = \alpha |z|^\ell$, where $\alpha \in \mathbb{C}$ is a constant. Consider the valuation

$$\varphi(K) = \int_{\partial K} \langle s, n(s) \rangle |s|^\ell d\sigma_K(s) = \beta \cdot \int_K |s|^\ell ds ,$$

where $\beta \in \mathbb{R} \backslash \{0\}$. Clearly, $(D\varphi)(K)(z) = \beta |z|^\ell \mathrm{vol}_2 K$. So $\Phi_2 = \frac{\alpha}{\beta} D\varphi$.    $\square$

We are going to describe O(2)-invariant polynomial valuations on $\mathcal{K}^2$ and SO($d$)-invariant polynomial valuations on $\mathcal{K}^d$, $d > 2$ (in the last case all of them turn out to be O($d$)-invariant).

Fix nonnegative integers $p$ and $q$. Consider for $K \in \mathcal{K}^d$

$$(4.7) \qquad \xi_{p,q}(K) = \int_{\partial K} \langle s, n(s) \rangle^p |s|^{2q} d\sigma_K(s) \ .$$

PROPOSITION 4.6.    *The function $\xi_{p,q}$ is an $\mathrm{O}(d)$-invariant continuous valuation of degree of polynomiality $\ell = p + 2q$ if $p \neq 1$, and of degree $\ell = 2q$ if $p = 1$.*

*Proof.* This is similar to the proof of Proposition 4.3.

*Remark.* Up to normalization $\xi_{1,q}$ coincides with $\int_K |s|^{2q} ds$.

Whenever we have the valuations $\xi_{p,q}$, we can consider "mixed" valuations

$$(4.8) \qquad \xi_{p,q}^{(j)}(K) = \frac{d^j}{d\varepsilon^j}\Big|_{\varepsilon=0} \xi_{p,q}(K + \varepsilon B) \ ,$$

where $B$ denotes the Euclidean ball, so that $\xi_{p,q}^{(j)}$ is also an $\mathrm{O}(d)$-invariant continuous polynomial valuation (here we again use Theorem 2.3).

THEOREM 4.7.    (a) *Every $\mathrm{O}(2)$-invariant continuous polynomial valuation on $\mathcal{K}^2$ is a linear combination of the $\xi_{p,q}^{(j)}$.*

(b) *If $d \geq 3$, then every $\mathrm{SO}(d)$-invariant continuous polynomial valuation on $\mathcal{K}^d$ is a linear combination of the $\xi_{p,q}^{(j)}$.*

*Remark.* It is easy to see that Theorem 4.7 is equivalent to Theorem B (i) in the introduction.

The next lemma will be needed in what follows.

LEMMA 4.8.    *Fix $\ell \geq 1$. Let $p \geq 1, q \geq 0$ be such that the valuation $\xi_{p,q}$ is polynomial of degree $\ell$ (see Proposition 4.6). Consider $\xi_{p,q}^{(j)}$ for $0 \leq j \leq d-2$ if $p > 1$, and $0 \leq j \leq d$ if $p = 1$. Then its image in $\Gamma_{d,\ell}$ (resp. in $\Gamma'_{d,\ell}$ if $d = 2$) can be described as follows:*

$$(4.9) \qquad \left(D\xi_{p,q}^{(j)}\right)(K)(x) = \binom{d-1}{j} |x|^{2q} \int_{S^{d-1}} \langle x, \omega \rangle^p dS_{d-1-j}(K, \omega) \quad \text{if} \quad p > 1 \ ,$$

*and*

$$(4.10) \qquad \left(D\xi_{1,q}^{(j)}\right)(K)(x) = \kappa W_j(K) |x|^{2q} \ ,$$

*where $W_j(K)$ is the $j^{\text{th}}$ quermassintegral and $\kappa > 0$ is a normalizing constant depending on $p, q, d, j$. Furthermore, all the $D\xi_{p,q}^{(j)}$ are linearly independent in $\Gamma_{d,\ell}$ for $p, q, j$ as above.*

*Proof.* As in the proof of Proposition 4.3 we can easily see that (4.9) and (4.10) hold for $j = 0$. Replacing $K$ by $K + \varepsilon B$ and taking derivatives with respect to $\varepsilon$, we obtain the general case.

Let us prove the linear independence. If some linear combination of valuations of types (4.9), (4.10) is zero, then we may assume that all the valuations included have the same degree of homogeneity with respect to $K$, say, $\mu$. Thus, for some $a_n$, $b$,

$$\sum_n a_n |x|^{2q_n} \int_{S^{d-1}} \langle x, \omega \rangle^{p_n} dS_\mu(K, \omega) + b|x|^\ell W_{d-\mu}(K) \equiv 0 \;,$$

where $2q_n + p_n = \ell$, $p_n > 1$ (note that in the case of odd $\ell$ the last summand disappears). Hence

$$\int_{S^{d-1}} \left( \left( \sum_n a_n |x|^{2q_n} \langle x, \omega \rangle^{p_n} \right) + b'|x|^\ell \right) dS_\mu(K, \omega) \equiv 0.$$

By an extension of Aleksandrov's theorem due to W. Weil ([W]), if $1 \leq \mu \leq d-1$ and a continuous function $f$ on $S^{d-1}$ satisfies $\int_{S^{d-1}} f(\omega) dS_\mu(K, \omega) = 0$ for all $K \in \mathcal{K}^d$, then $f$ has the form $f(\omega) = \langle a, \omega \rangle$ for some $a \in \mathbb{C}^d$ (Aleksandrov showed this for $\mu = d-1$).

Therefore, in our case, there exists $a(x) \in \mathbb{C}^d$, such that

$$b'|x|^\ell + \sum_n a_n |x|^{2q_n} \langle x, \omega \rangle^{p_n} \equiv \langle a(x), \omega \rangle \;.$$

Recall that for all $n$, the $p_n > 1$ are different. This implies that $b' = a_n = 0$.  □

As before, Theorem 4.7 follows by induction in $\ell$ from the following:

PROPOSITION 4.9.    *Let $p, q, \ell, j$ be as in Lemma 4.8.*

(a) *If $d = 2$, $\Gamma'_{2,\ell}$ is spanned by the $D\xi^{(j)}_{p,q}$;*

(b) *If $d \geq 3$, $\Gamma_{3,\ell}$ is spanned by the $D\xi^{(j)}_{p,q}$.*

*Proof.* By Theorem 2.5, it is sufficient to consider valuations of a given degree of homogeneity $\mu$, with $0 \leq \mu \leq d$.

*Case 1.*   $\mu = 0$. By Theorem 2.6(a), every such valuation $\Phi$ satisfies $\Phi(K) \equiv A$, $A \in T_{d,\ell}$, $\pi(U)A = A$ for all $U \in \mathrm{O}(2)$ if $d = 2$, or $U \in \mathrm{SO}(d)$ if $d \geq 3$. If $\ell$ is odd, then $A$ must be identically 0. But if $\ell$ is even, then $A$ is a polynomial, proportional to $|x|^\ell$. Thus, by Lemma 4.8, $\Phi(K) \equiv \left( D\xi^{(d)}_{1,\frac{\ell}{2}} \right)(K)$.

*Case 2.*   $\mu = d$. By Theorem 2.6(b), if $\Phi$ is homogeneous of degree $d$, then $\Phi(K) = B \cdot \mathrm{vol}_d(K)$. As in the previous case, if $\ell$ is odd, $B$ must be 0, and if $\ell$ is even, $B$ is proportional to $|x|^\ell$. But

$$|x|^\ell \cdot \mathrm{vol}_d(K) = \left(D\xi_{1,\frac{\ell}{2}}^{(0)}\right)(K) \ .$$

*Case* 3. $\mu = d - 1$. Using Theorem 2.6(c) and Remark 1 after it, we can see that there is one-to-one correspondence between valuations from $\Gamma'_{2,\ell}$ if $d = 2$ or $\Gamma_{d,\ell}$ if $d > 2$ and continuous functions

$$F : S^{d-1} \longrightarrow T_{d,\ell}$$

which satisfy

(i) $F$ is orthogonal to every linear functional on $S^{d-1}$,

(ii) $F(U\omega) = (\pi(U)F)(\omega)$ for all $\omega \in S^{d-1}$ and $U \in \mathrm{O}(2)$ if $d = 2$ or $U \in \mathrm{SO}(d)$ if $d > 2$. Again, we consider an intertwining operator

$$\widetilde{F} : T_{d,\ell}^* \longrightarrow C(S^{d-1})$$

defined as $\widetilde{F}(\xi)(\omega) = \langle \xi, F(\omega) \rangle$ for every $\xi \in T_{d,\ell}^*$, $\omega \in S^{d-1}$.

As before, $T_{d,\ell}^* = \sigma_d^\ell \oplus \sigma_d^{\ell-2} \oplus \sigma_d^{\ell-4} \oplus \ldots$, where all $\sigma_d^j$ are irreducible and pairwise nonequivalent. By Schur's lemma, $\widetilde{F}(\sigma_{\underline{d}}^j) \subset \sigma_d^j$. The condition (i) of orthogonality is clearly equivalent to saying that $\widetilde{F}(\sigma_d^1) = \{0\}$, which is satisfied automatically if $\ell$ is even (since $\sigma_d^1$ is not included in the decomposition of $T_{d,\ell}^*$).

By Schur's lemma the dimension of the linear space of all such intertwining operators is equal to

$$\frac{\ell}{2} + 1 \quad \text{if} \quad \ell \quad \text{is even, and}$$

$$\frac{\ell - 1}{2} \quad \text{if} \quad \ell \quad \text{is odd} \ .$$

Let us compute the dimension of the valuations spanned by the $D\xi_{p,q}^{(j)}$. First assume that $\ell$ is even. By Lemma 4.8, we have

$$(D\xi_{p,q}^{(0)})(K)(x) = \binom{d-1}{j}|x|^{2q} \int\limits_{S^{d-1}} \langle x, \omega \rangle^p dS_{d-1}(K, \omega)$$

if $p > 1$, $p + 2q = \ell$, and

$$\left(D\xi_{1,\frac{\ell}{2}}^{(1)}\right)(K)(x) = \kappa W_1(K)|x|^\ell \ ,$$

and these valuations are linearly independent. Therefore, the dimension of the linear span is equal to $\frac{\ell}{2} + 1$. Now assume $\ell$ to be odd. Again by Lemma 4.8 we have $\frac{\ell-1}{2}$ linearly independent valuations

$$\left(D\xi_{p,q}^{(0)}\right)(K)(x) = \binom{d-1}{j}|x|^{2q} \int\limits_{S^{d-1}} \langle x, \omega \rangle^p dS_{d-1}(K, \omega) \ ,$$

where $p > 1$, $p + 2q = \ell$. Thus this implies Case 3 and hence Proposition 4.9(a).

*Case* 4.    $\mu = 1$. Since for $d = 2$ the proposition follows from the previous cases, assume that $d > 2$. Fix $\Phi \in \Gamma_{d,\ell}$ such that $\Phi$ is homogeneous of degree 1. It is well-known ([H2]) that $\Phi$ must be Minkowski additive, i.e. $\Phi(\lambda_1 K_1 + \lambda_2 K_2) = \lambda_1 \Phi(K_1) + \lambda_2 \Phi(K_2)$ for all $K_i \in \mathcal{K}^d$, $\lambda_i \geq 0$. Since every $C^\infty$-function on $S^{d-1}$ is a difference of two smooth supporting functionals of two convex sets ([A]), $\Phi$ can be extended by linearity to a map

$$\Phi : C^\infty(S^{d-1}) \longrightarrow T_{d,\ell} ,$$

which clearly will be continuous (indeed, if $f_n \to f$ in the $C^\infty$-topology, then one can choose a large constant $M$ such that, for all $n$, the functions $f_n + M$ and $f + M$ will be supporting functionals of convex bodies). Moreover, $\Phi$ must be an intertwining operator of the quasi-regular representation $\pi$ of $\mathrm{SO}(d)$. Since

$$C^\infty(S^{d-1}) = \sigma_d^0 \oplus \sigma_d^1 \oplus \cdots \oplus \sigma_d^k \oplus \ldots ,$$
$$T_{d,\ell} = \sigma_d^\ell \oplus \sigma_d^{\ell-2} \oplus \ldots ,$$

Schur's lemma again implies that $\Phi(\sigma_d^k) \subset \sigma_d^k$. Translation invariance of $\Phi$ is equivalent to the property

$$\Phi(\sigma_d^1) = \{0\} ,$$

which is automatically satisfied for even $\ell$.

By Schur's lemma, $\Phi$ must be a composition of the orthogonal projection from $C^\infty(S^{d-1})$ onto $\sigma_d^\ell \oplus \sigma_d^{\ell-2} \oplus \ldots$ and some finite-dimensional operator, which is a multiplication by scalar operator on each component $\sigma_d^\ell, \sigma_d^{\ell-2}, \ldots$ and 0 on $\sigma_d^1$. Then obviously $\Phi$ has a unique continuous extension to the operator

$$\Phi : C(S^{d-1}) \longrightarrow T_{d,\ell} ,$$

which has the above properties. Moreover, every such operator restricted to the cone of supporting functionals of convex sets provides an example of continuous valuation from $\Gamma_{d,\ell}$, homogeneous of degree 1. If $\ell$ is even then the dimension of the linear space of such operators is equal to $\frac{\ell}{2} + 1$, and if $\ell$ is odd, to $\frac{\ell-1}{2}$. Similarly to Case 3, we see that the span of $\left\{ D\xi_{p,q}^{(d-2)} \right\}$ has the same dimension, where $p > 1$, $p + 2q = \ell$ if $\ell$ is odd, and the span of $D\xi_{1,\frac{\ell}{2}}^{(d-1)}$ and $\left\{ D\xi_{p,q}^{(d-2)} \right\}$ with $p > 1$, $p + 2q = \ell$ if $\ell$ is even.

Thus Case 4 is proved. This implies the proposition in the three-dimensional case.

*Case* 5.    Now let us assume that $d \geq 4$. It remains to consider the valuations from $\Gamma_{d,\ell}$ of degree of homogeneity $\mu$, $2 \leq \mu \leq d - 2$. Using $\{D\xi_{p,q}^{(j)}\}$ for appropriate $p, q, j$ as before, we see that the dimension of the linear space of these valuations is at least $\frac{\ell}{2} + 1$ if $\ell$ is even, and at least $\frac{\ell-1}{2}$ if $\ell$ is odd. We will show by induction in $d$ that these numbers also provide an upper estimate on the dimension; this will complete the proof of the proposition (the base $d = 3$ of induction is proved).

Let us fix some orthogonal decomposition $\mathbb{R}^d = \mathbb{R}^{d-1} \oplus \mathbb{R}^1$. Then

$$(4.11) \qquad T_{d,\ell} = T_{d-1,\ell} \oplus x_d \cdot T_{d-1,\ell-1} \oplus \cdots \oplus x_d^\ell \cdot T_{d-1,0} .$$

Consider the linear map

$$N : \Gamma_{d,\ell} \longrightarrow \Gamma_{d-1,\ell}$$

defined as follows: For every $\Phi \in \Gamma_{d,\ell}$ and for every compact convex subset $K \subset \mathbb{R}^{d-1}$, let

$$(N\Phi)(K) := \mathrm{Pr}_{d-1,\ell}(\Phi(K)) ,$$

where $\mathrm{Pr}_{d-1,\ell}$ is a projection from $T_{d,\ell}$ onto $T_{d-1,\ell}$ vanishing on the other summands of the decomposition $(4.11)$ (i.e. it is just a restriction to $x_d = 0$). Using the inductive assumption it is sufficient to show that $N$ is injective. Therefore suppose that $\Phi \in \mathrm{Ker} N$ and that $\Phi$ is homogeneous of degree $\mu$, for some $2 \le \mu \le d - 2$. Denote $S^{d-2} = S^{d-1} \cap \mathbb{R}^{d-1}$. Then for every $K \subset \mathbb{R}^{d-1}$ such that $\dim K \le d - 2$ we have $\Phi(K)\big|_{S^{d-2}} = 0$. Hence, by invariance of $K$ with respect to rotations about $(\mathrm{aff} \ K)^\perp$ and rotation equivariance of $\Phi$, we obtain that $\Phi(K) = 0$.

We will show that $\Phi$ is a simple valuation and this and Theorem 2.7 will conclude the proof. Let $\mathcal{K}^{d-1}$ denote the family of all compact convex subsets of $\mathbb{R}^{d-1}$. The restriction of $\Phi$ to $\mathcal{K}^{d-1}$ is simple; hence it suffices to check that $\Phi$ vanishes on orthogonal simplices in $\mathbb{R}^{d-1}$. Note that Theorem 2.7 and our assumption imply that the restriction of $\Phi$ to $\mathcal{K}^{d-1}$ (and hence $\Phi$ itself) is homogeneous of degree $\mu = d - 2$.

For every orthogonal simplex $S \subset \mathbb{R}^{d-1}$ we have a canonical decomposition ([H2]) of the simplex homothetic to $S$ with coefficient 2:

$$2 \cdot S = \bigcup_{j=0}^{d-1} (S_j' + S_{d-1-j}'') ,$$

where $S_j'$ and $S_{d-1-j}''$ are $j$- and $(d-1-j)$-dimensional orthogonal simplices in correspondence, lying in pairwise orthogonal subspaces. Thus

$$2^\mu \Phi(S) = \Phi(S_{d-1}') + \Phi(S_{d-1}'') + \sum_{j=1}^{d-2} \Phi(S_j' + S_{d-1-j}'')$$

$$= 2\Phi(S) + \sum_{j=1}^{d-2} \Phi(S_j' + S_{d-1-j}'') .$$

Since $\mu = d - 2 > 1$, it is sufficient to show that $\Phi(S'_j + S''_{d-1-j}) = 0$ (note that this trick was used in [Sch2]).

Now let us fix an orthogonal decomposition $\mathbb{R}^{d-1} = E_1 \oplus E_2$. Let $\dim E_i = k_i > 0$. We will choose coordinates in $\mathbb{R}^{d-1}$ such that $E_1 = \{(x_1, \ldots, x_{k_1})\}$, $E_2 = \{(x_{k_1+1}, \ldots, x_{d-1})\}$. Let $E'_2 = \{(x_{k_1+1}, \ldots, x_{d-1}, x_d)\}$. Thus $k'_2 := \dim E'_2 = k_2 + 1$.

By Theorem 2.3, for every $K_1 \subset E_1$, $K_2 \subset E_2$ we have a polynomial expansion, homogeneous of degree $d - 2$:

$$\Phi(\lambda_1 K_1 + \lambda_2 K_2) = \sum_{j=1}^{d-2} \lambda_1^j \lambda_2^{d-2-j} \Psi_j(K_1, K_2) \ ,$$

where $\lambda_i \geq 0$. Clearly, the coefficients $\Psi_j$ are simple translation invariant valuations in $K_1 \subset E_1$ and $K_2 \subset E_2$ separately. Moreover, $\Psi_j$ is homogeneous of degree $j$ in $K_1$ and of degree $d - 2 - j$ in $K_2$. It follows from Theorem 2.7 that if $\Psi_j \not\equiv 0$, then either $j = k_1$ and $d - 2 - j = k_2 - 1$ or $j = k_1 - 1$ and $d - 2 - j = k_2$. Hence

$$\Phi(\lambda_1 K_1 + \lambda_2 K_2) = \lambda_1^{k_1} \lambda_2^{k_2-1} \Psi_{k_1}(K_1, K_2) + \lambda_1^{k_1-1} \lambda_2^{k_2} \Psi_{k_1-1}(K_1, K_2) \ .$$

Because of the symmetry in $K_1$ and $K_2$, it suffices to prove that $\Psi_{k_1} \equiv 0$. If $\Psi_{k_1} \not\equiv 0$, then $k_2 - 1 > 0$. Now let $K_1 \subset E_1$ and $K'_2 \subset E'_2$. Then for $\lambda_1, \lambda_2 \geq 0$

$$\Phi(\lambda_2 K_1 + \lambda_2 K'_2) = \sum_{j=0}^{d-2} \lambda_1^j \lambda_2^{d-2-j} \Psi'_j(K_1, K'_2) \ .$$

Clearly if $K'_2 \subset E_2$, then $\Psi_{k_1}(K_1, K'_2) = \Psi'_{k_1}(K_1, K'_2)$. Let us fix $K_1 \subset E_1$ and consider $\Psi'_{k_1}(K_1, K'_2)$ as a valuation in $K'_2 \subset E'_2$. To shorten the notation let us denote $\Psi(K) := \Psi'_{k_1}(K_1, K)$. Clearly, $\Psi$ is an $SO(k'_2)$-equivariant valuation that is homogeneous of degree $k'_2 - 2$ which takes values in

$$(4.12) \qquad\qquad \bigoplus_{|\alpha| \leq \ell} x_1^{\alpha_1} \cdots x_{k_1}^{\alpha_{k_1}} \cdot T_{k'_2, \ell - |\alpha|} \ ,$$

where the sum extends over all the multi-indices $\alpha = (\alpha_1, \ldots, \alpha_{k_1})$ such that $|\alpha| := \sum_1^{k_1} \alpha_i \leq \ell$, and $T_{k'_2, m}$ denotes the space of homogeneous polynomials of degree $m$ in $k'_2 = k_2 + 1$ variables $x_{k_1+1}, \ldots, x_{d-1}, x_d$. Fix such an $\alpha$. Let $\Psi_\alpha$ be a projection of $\Psi$ on $T_{k'_2, \ell - |\alpha|}$ using (4.12). Thus $\Psi = \sum_{|\alpha| \leq \ell} x_1^{\alpha_1} \cdots x_{k_1}^{\alpha_{k_1}} \cdot \Psi_\alpha$. Since

$$(4.13) \quad T_{k'_2, \ell - |\alpha|} = T_{k'_2-1, \ell - |\alpha|} \oplus x_d \cdot T_{k'_2-1, \ell-|\alpha|-1} \oplus \cdots \oplus x_d^{\ell-|\alpha|} \cdot T_{k'_2-1, 0} \ ,$$

our assumption that $\Phi \in \mathrm{Ker} N$ implies that the projection of $\Psi_\alpha(K)$ on $T_{k'_2-1, \ell-|\alpha|}$ vanishes whenever $K \subset E_2$. But $\Psi_\alpha$ is homogeneous of degree

$k_2' - 2$ and $k_2' = k_2 + 1 \geq 3$. By the inductive hypothesis $\Psi_\alpha$ has the form

$$\Psi_\alpha(K) = \int\limits_{S^{k_2'-1}} F(\omega) dS_{k_2'-2}^{(k_2')}(K, \omega) \ ,$$

where $F : S^{k_2'-1} \to T_{k_2', \ell - |\alpha|}$ is a continuous function such that $F$ is orthogonal to every linear functional on the $S^{k_2'-1}$, and for every $U \in \mathrm{SO}(k_2')$ and every $\omega \in S^{k_2'-1}$,

$$F(U\omega) = (\pi(U)F)(\omega) \ .$$

Our proposition is immediate from the following:

LEMMA 4.10. *Let $k \geq 3$ and $F : S^{k-1} \to T_{k,m}$ be a continuous function. Fix a $(k-1)$-dimensional subspace $E := \{(y_1, \ldots, y_{k-1}, 0)\}$ and consider the decomposition*

$$T_{k,m} = T_{k-1,m} \oplus y_k \cdot T_{k-1,m-1} \oplus \cdots \oplus y_k^m \cdot T_{k-1,0} \ .$$

*Assume that*

(i) *$F$ is orthogonal to every linear functional on the sphere $S^{k-1}$;*

(ii) *$F(U\omega) = (\pi(U)F)(\omega)$ for every $U \in \mathrm{SO}(k), \omega \in S^{k-1}$;*

(iii) *for every convex compact set $K \subset E$ the projection*

$$\mathrm{Pr}_{k-1,m}\{\int\limits_{S^{k-1}} F(\omega) dS_{k-2}(K, \omega)\}$$

*vanishes (here $\mathrm{Pr}_{k-1,m}$ denotes the projection from $T_{k,m}$ to $T_{k-1,m}$ by the above decomposition).*

*Then $F \equiv 0$.*

*Proof.* First we observe that if $k \geq 3$, then (ii) must be true for every $U \in \mathrm{O}(k)$. Let $\omega_0$ be the point $(0, \ldots, 0, 1)$. One can easily check (e.g. using approximation by polytopes) that for every continuous function $f$ on the sphere $S^{k-1}$ and every $K \subset E (\simeq \mathbb{R}^{k-1})$

$$\int\limits_{S^{k-1}} f(\omega) dS_{k-2}(K, \omega) = \kappa \int\limits_{S^{k-2}} \left[ \int\limits_{[\omega_0, u, -\omega_0]} f(\omega) d\ell(\omega) \right] dS_{k-2}^{(k-1)}(K, u) \ ,$$

where $[\omega_0, u, -\omega_0]$ denotes the unique geodesic semicircle connecting $\omega_0$ and $-\omega_0$ and passing through $u \in S^{k-2} := S^{k-1} \cap E$. The inner integral in the right-hand side is taken with respect to the standard Lebesgue measure on this semicircle. The $S_{k-2}^{(k-1)}(K, u)$ denotes the $(k-1)$-dimensional surface area measure of $K$ on $S^{k-2}$; $\kappa$ is a normalization constant.

Thus for every $K \subset E$,

$$(4.14) \qquad \int\limits_{S^{k-2}} \mathrm{Pr}_{k-1,m} \left\{ \int\limits_{[\omega_0 u, -\omega_0]} F(\omega) d\ell(\omega) \right\} dS_{k-2}^{(k-1)}(K, u) = 0 \ .$$

Let us denote $G(u) := \mathrm{Pr}_{k-1,m} \left\{ \int\limits_{[\omega_0, u, -\omega_0]} F(\omega) d\ell(\omega) \right\}$. Then $G : S^{k-2} \to T_{k-1,m}$ is a continuous $\mathrm{O}(k-1)$-equivariant function. By (4.14)

$$\int\limits_{S^{k-2}} G(u) dS_{k-2}^{(k-1)}(K, u) = 0$$

for every $K \subset E$. Hence by Aleksandrov's theorem, for every $\xi \in T_{k-1,m}^*$ there exists a vector $a(\xi) \in E$ such that $\langle \xi, G(u) \rangle = (a(\xi), u)$. Hence, in fact, $G$ takes values in the subspace of $T_{k-1,m}$ of the spherical harmonics of degree one, which is just the space of $\mathbb{C}$-valued linear functionals on $E$, i.e. $(E \otimes \mathbb{C})^*$. Using the scalar product, we identify it with $E \otimes \mathbb{C}$. Since $G(u)$ is $\mathrm{O}(k-1)$-equivariant, it must have the form $G(u) = \alpha u$, where $\alpha \in \mathbb{C}$ is a fixed constant. Let $\widetilde{G}(\omega) := \alpha' \omega$, $\omega \in S^{k-1}$, where $\alpha'$ is a constant such that $\int\limits_{[\omega_0, u, -\omega_0]} \alpha' \omega d\ell(\omega)$ $= \alpha u$ (existence of $\alpha'$ is obvious). Setting $\widetilde{F}(\omega) := F(\omega) - \widetilde{G}(\omega)$, we see that for every $u \in S^{k-2}$,

$$(4.15) \qquad \mathrm{Pr}_{k-1,m} \int\limits_{[\omega_0, u, -\omega_0]} \widetilde{F}(\omega) d\ell(\omega) \equiv 0 \ .$$

Let us show that $\widetilde{F} \equiv 0$ (this and condition (i) of the lemma will finish the proof). Since $\widetilde{F}$ is $\mathrm{O}(k)$-equivariant, $\widetilde{F}(\omega_0)$ is invariant with respect to rotations about $\omega_0$. Hence it has the form

$$\widetilde{F}(\omega_0)(y_1, \ldots, y_{k-1}, y_k) = \sum_{j=0}^{[m/2]} c_j \left( \sum_{i=1}^{k-1} y_i^2 \right)^j y_k^{m-2j} \ ,$$

where $c_j \in \mathbb{C}$. For every $\theta \in [0, \pi]$, consider the point

$$\omega(\theta) = (0, \ldots, 0, -\sin\theta, \cos\theta) \ .$$

Clearly, $\omega(\theta)$ parametrizes the semicircle connecting $\omega_0$ and $-\omega_0$ and passing through the point $(0, \ldots, 0, -1, 0) \in E$. Hence

$\widetilde{F}(\omega(\theta))(y_1, \ldots, y_k)$

$$= \sum_{j=0}^{[m/2]} c_j \left( \sum_{i=1}^{k-2} y_i^2 + \left( y_{k-1}\cos\theta + y_k\sin\theta \right)^2 \right)^j (-y_{k-1}\sin\theta + y_k\cos\theta)^{m-2j} \ .$$

Condition (4.15) means that if we substitute $y_k = 0$, then

$$\int_0^\pi \sum_{j=0}^{[m/2]} c_j \left( \sum_{i=1}^{k-2} y_i^2 + y_{k-1}^2 \cos^2 \theta \right)^j (-y_{k-1} \sin \theta)^{m-2j} d\theta \equiv 0 \ .$$

Hence $c_j = 0$ for all $j$, i.e. $\widetilde{F} \equiv 0$.    $\square$

*Remark.* It follows from the proofs of Theorem 4.7 and Proposition 4.9 that the dimension of the space $\Omega'_{d,\ell}$ of all continuous O($d$)-invariant valuations, polynomial of degree at most $\ell$ (recall that if $d \geq 3$ then $\Omega'_{d,\ell} = \Omega_{d,\ell}$) satisfies the following recursive formula for $d \geq 2$:

$$\dim \Omega'_{d,\ell} = \dim \Omega'_{d,\ell-1} + \dim \Gamma'_{d,\ell} \ , \ \text{and}$$

$$\dim \Gamma'_{d,\ell} = (d-1) \cdot \left[ \frac{\ell}{2} \right] + \begin{cases} d+1 \ , & \text{if } \ell \text{ is even} \\ 0 \ , & \text{if } \ell \text{ is odd} \end{cases}$$

and $\dim \Omega'_{d,0} = d+1$ by the Hadwiger theorem.

Using the proofs of Theorem 4.4 and Lemma 4.5 similarly one gets a formula for the dimension of the space $\Omega_{2,\ell}$ of continuous SO(2)-invariant valuations on the plane, polynomial of degree at most $\ell$:

$$\dim \Omega_{2,\ell} = \dim \Omega_{2,\ell-1} + \begin{cases} \ell+3 \ , & \text{if } \ell \text{ is even} \\ \ell-1 \ , & \text{if } \ell \text{ is odd,} \end{cases} \ \text{and}$$

$$\dim \Omega_{2,0} = 3 \ .$$

Now we will consider general (nonpolynomial) continuous rotation invariant valuations and prove Theorem A of the introduction.

*Proof of Theorem* A. The case $d = 1$ follows from Proposition 3.1, and so let $d \geq 2$.

Let $C(\mathbb{R}^d)$ be the Fréchet space of complex-valued continuous functions on $\mathbb{R}^d$ with the usual topology of uniform convergence on compact sets.

Let $\mathcal{G}$ denote either SO($d$) or O($d$). For every continuous $\mathcal{G}$-invariant valuation $\varphi : \mathcal{K}^d \to \mathbb{C}$ consider the $C(\mathbb{R}^d)$-valued valuation $\Phi : \mathcal{K}^d \to C(\mathbb{R}^d)$ defined by $\Phi(K)(x) = \varphi(K+x)$. Then obviously $\Phi$ satisfies

(i) $\Phi$ is a continuous valuation;

(ii) $\Phi(K+y)(x) = \Phi(K)(x+y)$ for every $K \in \mathcal{K}^d$, $x, y \in \mathbb{R}^d$;

(iii) $\Phi(UK)(x) = \Phi(K)(U^{-1}x)$ for every $U \in \mathcal{G}$.

Since $\varphi(K) = \Phi(K)(0)$, we have a one-to-one correspondence between $\mathbb{C}$-valued continuous $\mathcal{G}$-invariant valuations and the space of $C(\mathbb{R}^d)$-valued valuations satisfying (i)–(iii). Denote the last space by $F$. Then obviously $F$ turns out to be the Fréchet space equipped with the following sequence of seminorms,

$$\|\Phi\|_N := \sup_{\substack{K \subset N \cdot B \\ |x| \leq N}} |\Phi(K)(x)| \, ,$$

for every $\Phi \in F$, and $N \in \mathbb{N}$, where $B$ is the Euclidean ball, and $|\cdot|$ is the Euclidean norm in $\mathbb{R}^d$. Note that $\{K \in \mathcal{K}^d \mid K \subset NB\}$ is a compact subset of $\mathcal{K}^d$ by the Blaschke selection theorem. Thus $\|\cdot\|_N$ is well-defined, i.e., finite. In $F$ we have a continuous representation $\pi$ of $\mathcal{G}$ given by

$$(\pi(U)\Phi)(K)(x) = \Phi(K)\Big(U^{-1}x\Big) \, .$$

It is well-known (see e.g. [He, Ch. IV, Lemma 1.9]) that for a continuous representation of any compact group $\mathcal{G}$ in Fréchet space the set of $\mathcal{G}$-finite vectors is dense in $F$ (recall that a vector $\xi$ is called $\mathcal{G}$-finite if span $(\mathcal{G}\xi)$ is finite dimensional). So to finish the proof it suffices to show that if $\Phi \in F$ is $\mathcal{G}$-finite valuation then $\Phi(K)(x)$ is a polynomial in $x$, whose degree is uniformly bounded in $K \in \mathcal{K}^d$. This follows from:

LEMMA 4.11. *Let $d \geq 2$. Let $E \subset C(\mathbb{R}^d)$ be a finite dimensional subspace, let $\dim E = \ell$ and assume that for every $a \in \mathbb{R}^d$, $E_a := \{f(x+a)|f(x) \in E\}$ is an* SO(d)*-invariant subspace. Then $E$ is contained in a subspace of polynomials of degree at most $d \cdot \ell$.*

*Proof.* We believe that this lemma is not new, but since we have no exact reference we present our proof for the convenience of the reader. First we reduce to the case of smooth functions. Let $\{\psi_n\}$ be sequence of SO(d)-invariant $C^\infty$-functions approximating the $\delta$-function at 0. Then $E_n = \{\psi_n * f | f \in E\}$ satisfies assumptions of the lemma and $E_n \subset C^\infty(\mathbb{R}^d)$. Now we may assume that $E \subset C^\infty(\mathbb{R}^d)$. Fix $a \in \mathbb{R}^d$, a skew-symmetric matrix $A$, and $f \in E$. Then for all $t \in \mathbb{R}$, $f\Big(e^{tA}x + a\Big) \in E_a$. Taking derivatives in $t$ at $t = 0$, we obtain

$$\langle df(x+a), Ax \rangle \in E_a \, .$$

Equivalently, $\langle df(x), Aa \rangle \in E$. But if $d \geq 2$, every vector $b \in \mathbb{R}^d$ can be represented as $b = Aa$ for some $a \in \mathbb{R}^d$ and a skew-symmetric matrix $A$. So $\langle df(x), b \rangle \in E$. In particular $\frac{\partial f}{\partial x_1} \in E$. Fix $f_1, \ldots, f_\ell$ a basis in $E$. Then for some matrix $B$ with constant entries

$$\begin{pmatrix} \frac{\partial f_1}{\partial x_1} \\ \cdot\cdot \\ \frac{\partial f_\ell}{\partial x_1} \end{pmatrix} = B \begin{pmatrix} f_1 \\ \cdot\cdot \\ f_\ell \end{pmatrix} \, .$$

This ordinary differential equation (where $x_2, \ldots, x_d$ are fixed) has a solution

$$f_i(x_1, x_2, \ldots, x_d) = \sum_{p,q=0}^{\ell-1} e^{\lambda_p x_1} x_1^q C_{p,q}(x_2, \ldots, x_d) \, ,$$

where the $\lambda_p \in \mathbb{C}$ are eigenvalues of $B$, and the $C_{p,q}$ are some functions of $x_2, \ldots, x_d$. However, $f_i$ has a similar expansion in $x_2, \ldots, x_d$. Finally, every $f \in E$ has the form

$$f(x) = \sum_i e^{\langle \xi_i, x \rangle} p_i(x) \ ,$$

where $p_i(x)$ are polynomials and $\xi_i \in \mathbb{R}^d \otimes \mathbb{C}$. Since $E$ is finite dimensional and $\mathrm{SO}(d)$-invariant, all $\xi_i$ must be equal to 0. $\qquad \square$

## 5. Applications to the integral-geometric formulas

Let us denote by $\mathcal{G}_{d,k}$ the Grassman manifold with the normalized Haar measure $\nu_k$ (thus $\nu_k$ is a probability measure on $\mathcal{G}_{d,k}$). Let $A_{d,k}$ be the manifold of $k$-dimensional *affine* subspaces of $\mathbb{R}^d$, and let $\mu_k$ be the Haar measure on $A_{d,k}$ with standard normalization (cf. [Sch1, p. 227]).

For every affine subspace $E \in A_{d,k}$ and every convex compact subset $M \subset E$, in this section we will denote by $M_\varepsilon$ the $\varepsilon$-extension of $M$ inside $E$. Clearly, $M_\varepsilon = M + \varepsilon \cdot (B \cap \overline{E})$, where $B$ is the Euclidean ball in $\mathbb{R}^d$, and $\overline{E}$ is the *linear* subspace, parallel to $E$. Then, by Theorem 2.3, $\int_{M_\varepsilon} |s|^2 dm_k(s)$ is a polynomial in $\varepsilon \geq 0$ of degree at most $k + 2$, where $m_k$ is the $k$-dimensional Lebesgue measure on $E$.

THEOREM 5.1.     *For every $K \in \mathcal{K}^d$, $1 \leq k \leq d - 1$,*

$$\int_{A_{d,k}} \frac{d^j}{d\varepsilon^j}\Big|_{\varepsilon=0} \Big[ \int_{(K \cap E)_\varepsilon} |s|^2 dm_k(s) \Big] d\mu_k(E)$$

$$= \begin{cases} \alpha_1 \cdot \frac{d^j}{d\varepsilon^j}\big|_{\varepsilon=0} \int_{K+\varepsilon B} |s|^2 dm_d(s) & \text{if} \quad j = 0, 1 \ ; \\[2mm] \alpha_2 \cdot \frac{d^j}{d\varepsilon^j}\big|_{\varepsilon=0} \int_{K+\varepsilon B} |s|^2 dm_d(s) + \beta_2 W_{j-2}(K) & \text{if} \quad 2 \leq j \leq k \ ; \\[2mm] \beta_3 W_{j-2}(K) & \text{if} \quad j = k+1, k+2 \ ; \\[2mm] 0 & \text{if} \quad j > k+2 \ ; \end{cases}$$

*where $\alpha_i$ and $\beta_i$ are constants depending on $d, k$ and $j$.*

*Proof.* The left-hand side in the theorem is easily seen to be a continuous $\mathrm{O}(d)$-invariant valuation in $K$, which is polynomial of degree at most 2 and homogeneous of degree $d + 2 - j$. Let us denote this valuation by $\varphi$ and compute its image $D\varphi$ in $\Gamma'_{d,2}$ (in the notation of Section 4). For every $K \in \mathcal{K}^d$,

$$\varphi(K + x) = \int_{A_{d,k}} \frac{d^j}{d\varepsilon^j}\big|_{\varepsilon=0} \Bigg[ \int_{(K \cap E)_\varepsilon} |s + x|^2 dm_k(s) \Bigg] d\mu_k(E)$$

$$= |x|^2 \int\limits_{A_{d,k}} \frac{d^j}{d\varepsilon^j}\Big|_{\varepsilon=0} \mathrm{vol}_k(K \cap E)_\varepsilon d\mu_k(E) + \text{lower order terms}$$

$$= c|x|^2 \int\limits_{A_{d,k}} W_j^{(k)}(K \cap E) d\mu_k(E) + \text{lower order terms} ,$$

where $W_j^{(k)}$ is the $k$-dimensional $j^{\text{th}}$ quermassintegral, and $c$ is a normalizing constant. The value of the last integral is well-known: it is the translation and rotation invariant continuous valuation in $K$, homogeneous of degree $d - j$, hence, by the Hadwiger theorem 2.4(c), it is proportional to $W_j(K)$. Thus

$$D\varphi(K) = c'|x|^2 W_j(K) \in \Gamma'_{d,2} .$$

By Lemma 4.8, the image of the valuation $\frac{d^j}{d\varepsilon^j}\Big|_{\varepsilon=0} \int\limits_{K+\varepsilon B} |s|^2 dm_d(s)$ is proportional to $D\varphi$. Hence $\varphi(K) - \alpha \cdot \frac{d^j}{d\varepsilon^j}\Big|_{\varepsilon=0} \int\limits_{K+\varepsilon B} |s|^2 dm_d(s)$ is a polynomial valuation of degree at most 1 for some constant $\alpha$. By Theorem 4.1, it must be translation invariant. Then it is proportional to $W_{j-2}$ by the Hadwiger theorem. □

Now let $E \in \mathcal{G}_{n,k}$ be a *linear $k$-dimensional subspace*. Denote $K|E$ the orthogonal projection of $K$ on $E$. As above, for every subset $M \subset E$ we will denote by $M_\varepsilon$ the $\varepsilon$-extension of $M$ *inside $E$*.

THEOREM 5.2. *For every $K \in \mathcal{K}^d$, $1 \leq k \leq d-1$,*

$$\int\limits_{\mathcal{G}_{d,k}} \frac{d^j}{d\varepsilon^j}\Big|_{\varepsilon=0} \left[ \int\limits_{(K|E)_\varepsilon} |s|^2 dm_k(s) \right] d\nu_k(E)$$

$$= \alpha \xi_{2,0}^{(d-1+j-k)}(K) + \beta \xi_{1,1}^{(d+j-k)}(K) + \gamma W_{d-2+j-k}(K) ,$$

*where $\xi_{\rho,q}^{(i)}$ are as in Section 4, and $\alpha, \beta, \gamma$ are constants depending on $d, k, j$.*

*Proof.* Again, the left-hand side is a continuous O($d$)-invariant valuation in $K$, polynomial of degree at most 2 and homogeneous of degree $k + 2 - j$. The rest is similar to the proof of the previous theorem. □

# 6. Several inequalities

THEOREM 6.1. *For every nonnegative integer $q$, the polynomial $\int\limits_{K+\varepsilon B} |s|^{2q} ds$ has nonnegative coefficients whenever $K \in \mathcal{K}^d$ and $K \ni 0$.*

*Proof.* We may assume that $K$ is a smooth, strictly convex body. Let $h_K$ be the support functional of $K$. One easily checks that

$$(6.1) \qquad \int_K |s|^{2q} ds = \frac{1}{d + 2q} \int_{\partial K} |s|^{2q} \langle s, n(s) \rangle d\sigma_K(s) \ ,$$

where $n(s)$ is the outer normal at $s \in \partial K$, $\sigma_K$ is the surface area measure on $\partial K$. Clearly, $\langle s, n(s) \rangle = h_K(n(s))$. Furthermore, it is well-known (c.f. e.g., [B, §94]) that the gradient of the support functional

$$\nabla h_K \Big|_{S^{d-1}} : S^{d-1} \longrightarrow \mathbb{R}^d$$

is the inverse of the Gauss-Bonnet map from $\partial K$ to $S^{d-1}$. Thus the right-hand side in (6.1) can be rewritten as

$$\frac{1}{d + 2q} \int_{S^{d-1}} |\nabla h_K(\omega)|^{2q} \langle \nabla h_K(\omega), \omega \rangle dS_{d-1}(K, \omega) \ .$$

Substituting $K + \varepsilon B$ for $K$ and using $\nabla h_{K + \varepsilon B}(\omega) = \nabla h_K(\omega) + \varepsilon \omega$ we obtain

$$\int_{K + \varepsilon B} |s|^{2q} ds$$

$$= \frac{1}{d + 2q} \int_{S^{d-1}} |\nabla h_K(\omega) + \varepsilon \omega|^{2q} \langle \nabla h_K(\omega) + \varepsilon \omega, \omega \rangle dS_{d-1}(K + \varepsilon B, \omega)$$

$$= \frac{1}{d + 2q} \int_{S^{d-1}} \left( |\nabla h_K(\omega)|^2 + \varepsilon^2 + 2\varepsilon \langle \nabla h_K(\omega), \omega \rangle \right)^q (\langle \nabla h_K(\omega), \omega \rangle + \varepsilon)$$

$$\cdot \sum_{i=0}^{d-1} \binom{d-1}{i} \varepsilon^i dS_{d-1-i}(K, \omega) \ .$$

Since $K \ni 0$, $\langle \nabla h_K(\omega), \omega \rangle = h_K(\omega) \geq 0$. Hence all the coefficients on the right-hand side of the last expression are nonnegative. $\qquad \square$

*Remark.* The valuation $\varphi(K) := \frac{d}{d\varepsilon} \Big|_{\varepsilon = 0} \int_{K + \varepsilon B} |s|^{2q} ds = \int_{\partial K} |s|^{2q} d\sigma_K(s)$ is "monotone" in the following sense: if $K_1 \supset K_2 \ni 0$, then $\varphi(K_1) \geq \varphi(K_2) \geq 0$. However, for higher derivatives in $\varepsilon$ $\frac{d^j}{d\varepsilon^j} \Big|_{\varepsilon = 0} \int_{K + \varepsilon B} |s|^{2q} ds$ (at least for even $j$) the similar monotonicity property on the class of convex compact sets containing 0 fails to be true (even in the 1-dimensional case). We do not know, whether this property holds on the class of centrally symmetric convex sets.

THEOREM 6.2. *Let* $K_1, \ldots, K_m \subset \mathbb{R}^2$ *be compact, convex,* centrally symmetric *subsets of the plane. Then* $\int_{\Sigma \lambda_i K_i} |s|^2 ds$ *is a polynomial in* $\lambda_i \geq 0$ *with* nonnegative *coefficients.*

*Remark.* We do not know what happens in higher dimensions.

*Proof of Theorem* 6.2. Clearly, the integral in the theorem is a homogeneous polynomial in $\lambda_i \geq 0$ of degree 4. Thus, analogously to the usual definition of mixed volumes, one can write

$$\int_{\Sigma \lambda_i K_i} |s|^2 ds = \sum_{i_1,i_2,i_3,i_4} P(K_{i_1}, K_{i_2}, K_{i_3}, K_{i_4}) \lambda_{i_1} \lambda_{i_2} \lambda_{i_3} \lambda_{i_4} ,$$

where $P(K_{i_1}, \ldots, K_{i_4})$ are coefficients not depending on any permutation of the indices $i_1, \ldots, i_4$. We have to show that $P(K_1, \ldots, K_4) \geq 0$ whenever the $K_i$ are centrally symmetric convex compact sets. As in the case of the mixed volumes, $P$ is Minkowski additive in each argument, i.e.

$$P(K_1' + K_1'', K_2, K_3, K_4) = P(K_1', \ldots, K_4) + P(K_1'', \ldots, K_4) .$$

It is known that every centrally symmetric convex compact subset of $\mathbb{R}^2$ is a zonoid (cf. [Sch1, Th. 3.5.1]). Thus it suffices to check the nonnegativity of $P$ for four centrally symmetric segments.

Fix $I_j = [-u_j, u_j]$, $j = 1, \ldots, 4$. We may choose the numeration of $u_j$ is such a way that all the vectors $u_1, \ldots, u_4$ lie in the same half-plane and we meet these vectors in this order if we move from the vector $u_1$ to $u_4$ counterclockwise (cf. Fig. 1).

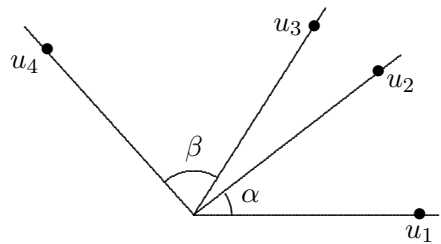

Fig. 1

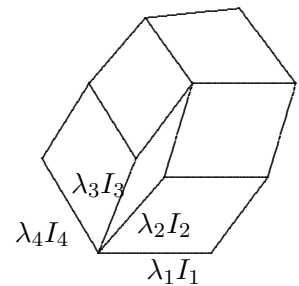

Fig. 2

Using the decomposition of the polygon $\sum_{j=1}^{4} \lambda_j I_j$ into parallelograms as drawn in Figure 2, one can easily see that

$$\int_{\sum_{1}^{4} \lambda_j I_j} |s|^2 ds = \int_{\lambda_1 I_1 + \lambda_2 I_2} |s - \lambda_3 u_3 - \lambda_4 u_4|^2 ds + \int_{\lambda_2 I_2 + \lambda_3 I_3} |s - \lambda_1 u_1 - \lambda_4 u_4|^2 ds$$

$$+ \int_{\lambda_1 I_1 + \lambda_3 I_3} |s + \lambda_2 u_2 + \lambda_4 u_4|^2 ds + \int_{\lambda_3 I_3 + \lambda_4 I_4} |s - \lambda_1 u_1 - \lambda_2 u_2|^2 ds$$

$$+ \int_{\lambda_2 I_2 + \lambda_4 I_4} |s + \lambda_3 u_3 - \lambda_1 u_1|^2 ds + \int_{\lambda_1 I_1 + \lambda_4 I_4} |s + \lambda_2 u_2 + \lambda_3 u_3|^2 ds \ .$$

A direct computation shows that

$$3 \cdot P(I_1, \ldots, I_4) = \langle u_3, u_4 \rangle u_1 \wedge u_2 + \langle u_1, u_4 \rangle u_2 \wedge u_3 - \langle u_2, u_4 \rangle u_1 \wedge u_3$$
$$+ \langle u_1, u_2 \rangle u_3 \wedge u_4 - \langle u_1, u_3 \rangle u_2 \wedge u_4 + \langle u_2, u_3 \rangle u_1 \wedge u_4 \ ,$$

where $\langle x, y \rangle$ denotes the scalar product of the vectors $x$ and $y$ and $x \wedge y$ denotes the oriented area of the parallelogram spanned by $x$ and $y$. It is easy to show that, for arbitrary $u_1, \ldots, u_4 \in \mathbb{R}^2$,

$$\langle u_1, u_4 \rangle u_2 \wedge u_3 - \langle u_2, u_4 \rangle u_1 \wedge u_3 - \langle u_1, u_3 \rangle u_2 \wedge u_4 + \langle u_2, u_3 \rangle u_1 \wedge u_4 = 0 \ .$$

Thus

$$3 \cdot P(I_1, \ldots, I_4) = \langle u_3, u_4 \rangle u_1 \wedge u_2 + \langle u_1, u_2 \rangle u_3 \wedge u_4 \ .$$

Let us denote by $\alpha$ the angle between $u_1$ and $u_2$, and by $\beta$ the angle between $u_3$ and $u_4$. By homogeneity, we may assume that $|u_1| = \cdots = |u_4| = 1$. Then the last expression can be rewritten

$$3 \cdot P(I_1, \ldots, I_4) = \cos \beta \sin \alpha + \cos \alpha \sin \beta = \sin(\alpha + \beta) \geq 0 \ . \qquad \square$$

## 7. Open questions

We would like to formulate a few questions, which are closely related to the material of this paper.

*Question* 1. Let $\varphi : \mathcal{K}^d \to \mathbb{R}$ be a valuation, continuous with respect to the Hausdorff metric and such that for every $K \in \mathcal{K}^d$, $\varphi(K+x)$ is a polynomial in $x$. Is it true that the degrees of these polynomials are *uniformly bounded* for all $K$?

(Note that in the 1-dimensional case the answer is positive, and it easily follows from the Bair theorem.)

*Question* 2 (due to V.D. Milman). Is it true that every continuous valuation can be approximated by polynomial valuations uniformly on compact subsets of $\mathcal{K}^d$ (without any assumption on rotation invariance)?

*Question* 3 (due to V.D. Milman). Let $\varphi(K) := \frac{d^j}{d\varepsilon^j}\Big|_{\varepsilon=0} \int_{K+\varepsilon B} |s|^{2q}ds$. Does the following monotonicity property hold on the class of centrally symmetric, convex, compact sets: if $K_1 \supset K_2$ are centrally symmetric, then $\varphi(K_1) \geq \varphi(K_2)$ (compare Theorem 6.1 and the remark following it)?

*Question* 4. Let $K_1, \ldots, K_s \subset \mathbb{R}^d$ be centrally symmetric compact convex sets. Let $q$ be a natural number. Under what condition has the polynomial $\int_{\sum_i \lambda_i K_i} |s|^{2q}ds$ in $\lambda_i \geq 0$ nonnegative coefficients? Is it true at least for $q = 1$ and $K_1, \ldots, K_s$ zonoids (compare Theorem 6.2)?

SCHOOL OF MATHEMATICAL SCIENCES, TEL AVIV UNIVERSITY, RAMAT-AVIV, ISRAEL
*E-mail address*: semyon@math.tau.ac.il